\documentclass[numbers,preprint,review,12pt]{elsarticle}

\usepackage{amssymb}
 \usepackage{amsthm}
\usepackage{amsgen,amsmath,amstext,amsopn,amsfonts,amsthm}
\journal{}

\begin{document}
\newcommand{\iref}[1]{\eqref{#1}}

\newtheorem{theorem}{Theorem}
\newtheorem{theoreme}{Theorem}
\newtheorem{remarque}{Remark}[section]
\newtheorem{lemme}{Lemma}[section]
\newtheorem{proposition}{Proposition}[section]
\newtheorem{definition}{Definition}[section]

\def \trait (#1) (#2) (#3){\vrule width #1pt height #2pt depth #3pt}

\def \fin{\hfill
        \trait (0.1) (5) (0)
        \trait (5) (0.1) (0)
        \kern-5pt
        \trait (5) (5) (-4.9)
        \trait (0.1) (5) (0)
\medskip}

\def\N{\Bbb N}
\def\R{\Bbb R}
\def\Po{\Bbb P}
\def\Z{\Bbb Z}
 
\def\C{\Bbb C}
\def\Q{\Bbb Q}
\def\T{\Bbb T}
\def\RR{\rm \hbox{I\kern-.2em\hbox{R}}}
\def\NN{\rm \hbox{I\kern-.2em\hbox{N}}}
\def\ZZ{\rm {{\rm Z}\kern-.28em{\rm Z}}}
\def\vp{\varphi}
\def\<{\langle}
\def\>{\rangle}
\def\t{\tilde}
\def\li{\infty}
\def\e{\varepsilon}
\def\lsim{\raisebox{-1ex}{$~\stackrel{\textstyle<}{\sim}~$}}
\def\gsim{\raisebox{-1ex}{$~\stackrel{\textstyle>}{\sim}~$}}
\def\sm{\setminus}
\def\nl{\newline}
\def\ov{\overline}
\def\wt{\widetilde}
\def\Chi{\raise .3ex \hbox{\large $\chi$}} 
\def\vp{\varphi}
\def\d{\delta}
\def\dt{\t\delta}
\def\a{\alpha}
\def\b{\beta}
\def\g{\gamma}
\def\v1{\tilde v}
\def\D{\Delta}
\def\phit{\tilde \varphi}
\def\[{\Bigl [}
\def\]{\Bigr ]}
\def\({\Bigl (}
\def\){\Bigr )}
\def\[{\Bigl [}
\def\]{\Bigr ]}
\def\({\Bigl (}
\def\){\Bigr )}
%%%%%%%%%%%%%%%%%%%%%%%%%%%%%%%%%%%%%%%%%%%%%%%%%%%%%%%%
\def\cM{{\cal M}}
\def\cC{{\cal C}}
\def\cO{{\cal O}}
%%%%%%%%%%%%%%%%%%%%%%%%%%%%%%%%%%%%%%%%%%%%%%%%%%%%%%%%%%
%%%notaions 
\def\bR{\mathbb{R}^{d}}
\def\bZ{\mathbb{Z}^{d}}
\def\bN{\mathbb{N}_0^{d}}
\def\lp{\ell^p(\mathbb{Z}^{d})}
\def\li{\ell^{\infty}(\mathbb{Z}^{d})}
\def\l0{\ell^0(\mathbb{Z}^{d})}
\def\Lp{L^p(\mathbb{R}^{d})}
\def\Li{L^{\infty}(\mathbb{R}^{d})}
\def\x{\bf x}
\def\k{\bf k}
%%%%%%%%%%%%%%%%%%%%%%%%%%%%%%%%%%%%%%%%%%%%%%%%%%%%%%%%%%%

\setcounter{page}{1}

\newcommand{\argm}{\rm argmin}
\newcommand{\supp}{\rm supp}

\begin{frontmatter}



\title{Smoothness of Nonlinear and Non-Separable Subdivision Schemes}
\author{Basarab Matei$^{(a)}$, Sylvain Meignen $^{*,(b)}$ and Anastasia
Zakharova$^{(b)}$}  
\address{\large 
\begin{tabular}{c}
\vspace{0.5 cm}\\
$^{(a)}$LAGA Laboratory, Paris XIII University, France,\\
Tel:0033-1-49-40-35-71\\
FAX:0033-4-48-26-35-68\\
E-mail: matei@math.univ-paris13.fr\\
\vspace{0.5 cm}\\
$^{(b)}$ LJK Laboratory, University of Grenoble, France\\ 
Tel:0033-4-76-51-43-95\\
FAX:0033-4-76-63-12-63\\
E-mail: sylvain.meignen@imag.fr\\
E-mail: anastasia.zakharova@imag.fr\\
\
\vspace {8 cm}
\end{tabular}\\
} 
\begin{abstract}
We study in this paper nonlinear subdivision schemes in a multivariate setting allowing 
arbitrary dilation matrix. We investigate the convergence of such iterative process to some limit function.
Our analysis is based on  some conditions on the contractivity of the associated scheme for the differences. 
In particular, we show the regularity of the limit function, in $L^p$ and Sobolev spaces.   
\end{abstract}

\begin{keyword}
Nonlinear subdivision scheme \sep  convergence of subdivision schemes \sep box splines\\

%% MSC codes hec    in \cite{HanJia} the convergence in $L^p$ is carried out and then 
%generalized in Sobolev spaces in \cite{Jia1} and \cite{Jia2}
% in the form: \MSC code \sep code

\end{keyword}

\end{frontmatter}


\newpage 

\section{Introduction}
Subdivision schemes have been the subject of active research in recent years.
In such algorithms, discrete data are recursively generated from coarse to fine by
means of local rules. When the local rules are independent of the data, the underlying 
refinement process is linear. This case is extensively studied in literature. 
The convergence of this process and the existence of the limit function  
was studied in \cite{CDM} and \cite{Dyn} when the 
scales are dyadic. When the scales are related to a dilation
matrix $M$, the convergence to a limit function in $L^p$ was studied in 
\cite{HanJia} and 
generalized to Sobolev spaces in \cite{Jia1} and \cite{Jia2}.
In the linear case, the stability is a consequence of 
 the smoothness  of the limit function.

The nonlinearity arises naturally when one needs to adapt locally the refinement 
rules to the data such as in image or geometry processing. 
Nonlinear subdivision schemes based
on dyadic scales were originally introduced by Harten 
\cite{harten}\cite{hartenENO} through the so-called essentially non-oscillatory (ENO) 
methods. These methods have recently been  
adapted to image processing into essentially non-oscillatory 
edge adapted (ENO-EA) methods. Different versions of ENO methods exist 
either based on
polynomial interpolation as in \cite{Co}\cite{Ar} or in a wavelet 
framework \cite{Chan}, corresponding to interpolatory or non-interpolatory 
subdivision schemes respectively.    

In the present paper, we study nonlinear subdivision schemes
associated to dilation matrix $ M$. After recalling the definitions 
on nonlinear subdivision schemes in that context, we give sufficient
conditions for convergence in Sobolev and $L^p$ spaces.  
\section{General Setting}
\subsection{Notations}
%%%%%%%%%%%%%%%%notations%%%%%%%%%%%%%%%%%% 

Before we start, let us introduce some notations that will be used throughout the paper. 
We denote $\# Q$ the cardinal of the set $Q$. 
For a multi-index $\mu=(\mu_1,\mu_2,\cdots,\mu_d) \in \mathbb{N}^{d}$ and a vector 
$x=(x_1,x_2,\cdots,x_d)\in \bR$   we define
$|\mu|=\sum\limits_{i=1}^d{\mu_i}$,  $\mu !=\prod\limits_{i=1}^d{\mu_i!}$ and  $x^{\mu}=\prod\limits_{i=1}^d{{ x_i }^{\mu_i}}.$

For two multi-index $m,\mu\in \mathbb{N}^{d}$ we also define 
$$\left (\begin{array}{c} 
\mu\\
m
\end{array}
\right )=\left (\begin{array}{c} 
\mu_1\\
m_1
\end{array}
\right )
\cdots 
\left ( \begin{array}{c} 
\mu_d\\
m_d
\end{array}
\right ).
$$
%%%%%%%%%%%%%%%%%esapces des suites %%%%%%%%%%%%%%%%%%%%%%%%%%%%%%%%%%%%%%%%%%%%%

Let  $\ell(\bZ)$ be the space of all sequences indexed by $\bZ$. 
The subspace of bounded sequences is denoted by 
$\li$ and  $\|u\|_{\li}$ is the supremum of $\{ |u_k|:k \in \bZ\}$. 
We denote $\l0$ the subspace of all sequences with finite 
support (i.e. the number of non-zero components of a sequence is finite).
As usual, let $\lp$ be
the Banach space of sequences
$u$ on $\bZ$ such that $\|u\|_{\lp} < \infty$, where
$$
\|u\|_{\lp} := \left ( \sum\limits_{k \in \bZ} |u_k|^p \right
)^{\frac{1}{p}} \textrm { for }1\leq p < \infty.
$$
%%%%%%%%%%%%%%%%esapces des fonctions%%%%%%%%%%%%%%%%%%%%%%%%%%%%

As in the discrete case, we denote by $\Lp$ the space of all
measurable functions $f$ such that $\|f\|_{\Lp} < \infty$, where
$$
\|f\|_{\Lp} := \left ( \int_{\bR} |f(x)|^p dx \right )^\frac{1}{p} 
\textrm{ for } 1 \leq p < \infty
$$
and $\|f\|_{\Li}$ is the essential supremum of $|f|$ on $\bR$.

%%%%%%%%%%%%%%%%%%%%%%%%%%%%%%%%%%%%%%%%%%%%%%%%%%%%%%%%%%%%%%

%%%%%%%%%%%%%%%%%Differnces et differential op 
Let $\mu\in \mathbb{N}^{d}$ be a multi-index, we define  $\nabla^{\mu}$ the
difference operator $\nabla_1^{\mu_1}\cdots \nabla_d^{\mu_d}$, where 
$\nabla_j^{\mu_j}$ is the $\mu_j$th 
difference operator with respect to the $j$th coordinate of the canonical basis. 
We define $D^{\mu}$ as  $D_1^{\mu_1}\cdots D_d^{\mu_d}$, where $D_j$ is the
differential operator with respect to the $j$th coordinate of the canonical basis. 
Similarly, for a vector $x \in \R^{d}$ the differential 
operator with respect to  $x$ is denoted by  $D_x$.
%%%%%%%%%%%%%%%%%%%%%%%%%%%%%%%%%%%%%%%%%%%%%%%%%%%%%%%%%%%%%%%%%%%%%%

%%%%%%%%%%%%matrice de dilations

A matrix $M$  is called a dilation matrix if it has integer entries and if  
$\lim\limits_{n\rightarrow \infty}M^{-n} =0$.
In the following, the invertible dilation matrix is always denoted by $M$ 
and $m$ stands for $|det(M)|$. 
%%%%%%%%%%%%%%%%%%%%%%%%%%%

%%%%%%%%%%%%%%%%fonctions %%%%%%%%%%%%%%%%%%%
For  a dilation matrix $M$  and any arbitrary function $\Phi$ we put 
$\Phi_{j,k}(x) = \Phi(M^{j}x -k)$. 
 
We also recall that  a compactly supported function 
$\Phi$ is  called $L^p$-stable if there exist two constants $C_1, C_2 >0$
satisfying  
$$
\quad C_1 \|c\|_{\lp} \leq \|\sum\limits_{k \in \Z^{d}} 
c_k \Phi(x - k)\|_{\Lp} \leq C_2 \|c\|_{\lp}.
$$
Finally, for two positive quantities $A$ and  $B$ depending on a set of parameters,
the relation $A\lsim B$  implies the  existence  of a  positive constant $C$,
independent of the parameters, such  that $A\leq C B$. Also $A\sim B$
means $A\lsim B$ and  $B\lsim A$.

\subsection{Local, Bounded and Data Dependent Subdivision Operators, Uniform Convergence
Definition}
In the sequel, we will consider the general class of local,  bounded and data dependent subdivision operators which are defined as follows:
\newtheorem{LocalPred}{Definition}
\begin{LocalPred}
\label{locpred}
For  $v \in \li$,   
a local, bounded and data dependent  subdivision operator is defined by 
\begin{eqnarray}\label{subdi}
S(v) w_k = 
\sum\limits_{l \in \bZ} a_{k-M l} (v) 
w_l,
\end{eqnarray}
for any $w$ in $\li$ and where the real coefficients $a_{k- M l}(v) \in \R$ 
are such that 
\begin{eqnarray}
\label{localdef}
a_{k- M l}(v) = 0,  \quad if  \quad \|k - M l \|_{\li} > K
\end{eqnarray}
for a fixed constant $K$.
The coefficients $a_{k}(v)$ are assumed to be 
uniformly bounded by a constant $C$, i.e. there is $C>0$ independent of $v$ 
such that:
$$
|a_k(v)| \leq C.
$$
Note that the definition of the coefficients  depends on some sequence 
$v$, while $S(v)$ acts on the sequence $w$.

Note also that, from  \iref{subdi} and \iref{localdef} 
the new defined value $S(v) w_k$ depends only on those values $l$ 
satisfying   $\|k - M l \|_{\li} > K$. The subdivision operator in this sense is  
{\it local}.
\end{LocalPred}

To simplify, in what follows a data dependent subdivision operator
is an operator in the sense of Definition \ref{locpred}.
With this definition, the associated subdivision scheme 
is the recursive action
of the data dependent  rule $Sv = S(v)v$ on an initial set of data $v^0$, 
according to:
\begin{eqnarray}
v^j = Sv^{j-1} = S(v^{j-1})v^{j-1}, \ j \geq 1.
\end{eqnarray}
\subsection{Polynomial Reproduction for Data Dependent Subdivision Operators}
The study of the convergence of data dependent subdivision operators 
will involve the polynomial reproduction property. We recall the definition of
the space $\mathbb{P}_N$ of polynomials of total 
degree $N$:
$$
\mathbb{P}_N:= \{ P; P(x) = \sum\limits_{|{\bf \mu}|\leq N}  a_\mu 
x^{\mu} \}.
$$
With these notations, the polynomial reproduction properties read: 
\newtheorem{Exactness}[LocalPred]{Definition}
\begin{Exactness}
Let $N \geq 0$ be a fixed integer. 
\begin{enumerate}
\item The data dependent subdivision operator $S$ has 
the property of reproduction of polynomials of total degree $N$ if for all $u \in 
\li$ and $P\in\mathbb{P}_N$ there exists $\tilde P \in \mathbb{P}_N$
 with $P-\tilde P \in \mathbb{P}_{N-1}$ such that $S(u)p=\tilde p$ where $p$
 and $\tilde p$ are defined by $p_k=P(k)$ and $\tilde p_k 
 = \tilde P (M^{-1}k)$.
\item The data dependent subdivision operator $S$ has 
the property of exact reproduction of polynomials of total degree $N$ if for all $u \in 
\li$ and $P\in\mathbb{P}_N$,  $S(u)p=\tilde p$ where $p$
 and $\tilde p$ are defined by $p_k=P(k)$ and $\tilde p_k 
 = P (M^{-1}k)$.
\end{enumerate}
\end{Exactness}
\underline{Remark:} 

The case $N=0$ is  the so-called "constant reproduction property". 
For a data dependent subdivision operator defined as in \iref{subdi}, 
the constant reproduction property  reads $\sum\limits_{k\in\bZ}  a_{k-M l}(v) = 1,$  for all  ${v} \in \li$.

\section{Definition of Schemes for the Differences}
Another ingredient for our study is the schemes for the differences associated
to the data dependent subdivision operator. 
The existence of schemes for the differences is obtained by using the 
polynomial reproduction property of the data dependent subdivision operator.

Let us denote  $\Delta^{l} = (\nabla^\mu,|\mu|=l)$ and then state 
the following result on the existence of schemes for the differences: 
 
\newtheorem{DataDependent}{Proposition}
\begin{DataDependent}
Let $S$ be a data dependent subdivision operator which reproduces polynomials up
to total degree $N$. Then for $1\leq l\leq N+1$ there exists a data dependent 
subdivision rule $S_l$ with the property that for all $v$,$w$ in $\li$, 
$$
\Delta^l S(v) w := S_l(v) \Delta^l w
$$
\end{DataDependent}
\textsc{Proof:} 
%%%%%%%%%%%%%%Nouvelle preuve

Let $l$ be an integer such that $1 \leq l \leq N+1$.
By using the definition of $\nabla^{\mu}$ with  $|\mu|=l$, we write:
\begin{eqnarray} 
\nonumber
\nabla^{\mu} S(v) w_k =
\nabla^{\mu_1}_1
\cdots
\nabla^{\mu_d}_d S(v) w_k. 
\end{eqnarray}
From the definition of $S(v) w$  we infer  that 
\begin{eqnarray} 
\nonumber
\nabla^{\mu} S(v) w_k = \sum\limits_{m_1,\cdots,m_d=0}^{\max (\mu_1,\cdots,\mu_d)}
(-1)^{l} \left 
(\begin{array}{c}
\mu\\
m
\end{array}
\right )
\sum_{p \in \mathbb{Z}^{d}}
a_{k-m\cdot e-Mp}(v)
w_p, 
\end{eqnarray}
where we have used the notation $m\cdot e= m_1e_1+\cdots+m_d e_d$.
Straightforward computations give
\begin{eqnarray} 
\nonumber
\label{vjj1}
\nabla^{\mu} S(v) w_k &=& \sum\limits_{p \in \mathbb{Z}^{d}}
w_p
\sum\limits_{m_1,\cdots,m_d=0}^{\max (\mu_1,\cdots,\mu_d)}
(-1)^{l} \left (\begin{array}{c}
\mu\\
m
\end{array}
\right )
a_{k - m\cdot e -Mp}(v)\ \\
&=&\sum\limits_{p \in \mathbb{Z}^{d}}
w_p f_{k,p}(v,\mu).
\end{eqnarray}
Let us clarify the definition of $f_{k,p}(v,\mu).$ 
Since  the data dependent subdivision operator is local we have $a_{k-Mp}(v) = 0$ for any data $v\in\li$  and any index 
$k$ such that $\|k-Mp\|_{\li} > K$. Now by putting  
$k = \varepsilon + Mn$, we get that $f_{k,p}(v,\mu)$ is defined for $p$ in the set
$$
V^{\mu}(k) :=
\left \{ p : \|n-p+{M}^{-1}(\varepsilon- m\cdot e)
\|_{\infty} \leq K \|M^{-1}\|_{\infty},\; 0 \leq m_i \leq \mu_i \forall i \right \}
$$
Then, we define
$V(k) := \left \{ p : \|k-Mp\|_{\infty} \leq K \right \} $.
Since the data dependent subdivision scheme reproduces polynomials up
to total degree $N$, we have for any $|\nu|= r \leq N$:
\begin{eqnarray}
\label{poly}
\sum\limits_{p \in V (k)}
a_{k-Mp}(v) p^{\nu} = P_{\nu}(k) \quad \mbox{ for all }k \in \Z^d,
\end{eqnarray}
where $P_{\nu}$ is a polynomial of total degree $r$.
By tacking the differences of order $|\nu'| = r+1$ in \iref{poly} we get
\begin{eqnarray*}
\sum\limits_{p \in V^{\nu'} (k)}
f_{k,p}(v,\nu') p^{\nu} &=& 0.
\end{eqnarray*}
Note that  the above equality is true for any $\nu$ 
such that $|\nu|=r$. We deduce  
that $\left ( f_{k,p}(v,\nu') \right )_{k} \in \Z^d$
is orthogonal to $\left ( p^q \right )_{p \in V^{\nu'}(k)}$  where 
$|q| \leq r$. Note that 
$
\left \{
\left ( \nabla^{\nu} \delta_{n-\beta}
\right )_{n \in V^{\nu'}(k)},
|\nu|= r+1,
\beta \in \Z^d
\right \}
$ spans $\left ( p^q \right )_{p \in V^{\nu'}(k)}$ 
and we may thus write for any $p \in V^{\nu'}(k)$:
\begin{eqnarray*}
f_{k,p}(v,\nu')=
\sum\limits_{|{\nu}| = r+1}
\sum\limits_{r \in \mathbb{Z}^{d}}
c^\mu_{k,r}(v) \nabla^{\mu} \delta_{p-r}.
\end{eqnarray*}
Now, by using  \iref{vjj1}  we obtain for any $|\mu| \leq N+1$:
\begin{eqnarray*}
\nabla^{\mu} S(v)w_k &=&
\sum\limits_{p \in V^{\mu}(k)}
w_p
\sum\limits_{|\nu| = l}
\sum\limits_{r \in \Z^d}
c^\nu_{k,r}(v) \nabla^{\nu} \delta_{p-r}\\
&=&\sum\limits_{ p \in V^{\mu}(k)}
\sum\limits_{|{\nu}| = l}
c^\nu_{k,r}(v)
\nabla^\nu w_p
\end{eqnarray*}
If we now make $\mu$ vary, we obtain the desired relation.
\fin \\
Now that we have proved the existence of schemes for the differences, 
we introduce the notion of joint spectral radius for these schemes, 
which is a generalization of the one dimensional case which 
can be found in \cite{Ost}.
\newtheorem{SpectralRad}[LocalPred]{Definition}
\begin{SpectralRad}
Let $S(v): \lp  \rightarrow \lp $ be a 
data dependent subdivision operator such that the difference operators 
$S_l(v): \left ( \lp \right )^{q_l} \rightarrow \left (\lp \right)^{q_l} $, 
with $q_l =\#\{\mu,|\mu|=l\}$ 
exists for $l \leq N+1$. 
Then, to each operator $S_l$, $l=0,\cdots,N+1$ (putting $S_0 = S$) 
we can associate  the joint spectral radius  given by 
$$
\rho_{p,l}(S) := \inf\limits_{j \geq 1} 
\|(S_l)^j\|_{\lp^{q_l}}^{\frac{1}{j}}. 
$$ 
\end{SpectralRad}
In other words, $\rho_p(S)$ is the infimum of all $\rho >0$ 
such that for all $v \in  \lp $, one has 
\begin{eqnarray}
\label{spectprat}
\|\Delta^l S^j v \|_{\lp^{q_l}} \lsim \rho^j \|\Delta^l v\|_{\lp^{q_l}},
\end{eqnarray}
for all $j \geq 0$.\\
\noindent \underline{Remark:}  Let us define a set of vectors 
$\{x_{\bf 1},\cdots,x_{\bf n}\}$ such that 
$[x_{\bf 1},\cdots,x_{\bf n}]\mathbb{Z}^n=\bZ$, 
$n \geq d$ (i.e. a set such that 
the linear combinations of its elements with coefficients in $\mathbb{Z}$ 
spans $\bZ$). We use the bold notation in the definition of the set so as 
to avoid the confusion with the coordinates of vector $x$. Then, consider the differences in the directions $x_{\bf 1},\cdots,x_{\bf n}$. One can show that there exists a scheme for that 
differences which we call $\tilde S_l$ for $l \leq N+1$ provided the data dependent 
subdivision operator reproduces polynomials up 
to degree $N$  (the proof is similar to that using the canonical directions). 
If we denote by $\tilde \Delta^l$ the difference operator of order $l$ in the
directions $x_{\bf 1},\cdots,x_{\bf n}$, one can see that 
$\|\tilde \Delta^l v\|_{\lp^{\tilde q_l}}  \sim  \|\Delta^l v\|_{\lp^{q_l}}$  
for all $v$ in $\lp$ and where 
$\tilde q_l = \# \{\mu, |\mu|=l, \mu =(\mu_i)_{i=1,\cdots,n} \}$.Then, 
following (\ref{spectprat}), one can deduce that the joint spectral radius 
of $\tilde S_l$ is the same as that of $S_l$. 
\section{Convergence in $L^p$ spaces}
In the following, we study the convergence of data dependent 
subdivision schemes in $L^p$ which corresponds to the following 
definition: 
\newtheorem{Converg}[LocalPred]{Definition}
\begin{Converg}
The subdivision scheme $v^j =Sv^{j-1}$ 
converges in $\Lp$, if for every set of initial control points 
$v^0 \in\lp$, there exists a non-trivial function $v$ in 
$\Lp$, called the limit function, such that 
$$
\lim\limits_{j \rightarrow \infty} \|v_j-v\|_{\Lp}=0.
$$
\end{Converg}
where $v_j(x) = \sum\limits_{k \in \bZ} 
v^j_k \phi_{j,k}(x)$ with 
$\phi(x) = \prod\limits_{i=1}^{d} \max(0,1-|x_i|)$.
\subsection{Convergence in the Linear Case}
When $S$ is independent of $v$, the rule (\ref{subdi}) defines a linear subdivision
scheme:
$$
Sv_k = 
\sum\limits_{l\in \bZ} a_{k-M l} 
v_l.
$$
If the linear subdivision scheme converges for any $v \in
\lp$ to some function in $\Lp$ and if there exists 
$v^0$ such that $\lim\limits_{j\rightarrow +\infty} v^j \neq 0$, 
then $\{a_k, k \in \bZ\}$ determines a unique continuous compactly supported
function $\Phi$ satisfying 
\begin{eqnarray*}
\Phi (x) = 
\sum\limits_{k \in \bZ} a_k \Phi(Mx-k) \textrm{ and }
\sum\limits_{k \in \bZ} \Phi (x -k) = 1.
\end{eqnarray*}
Moreover, $v(x)= \sum\limits_{k \in \bZ}
v^0_k \Phi (x-k)$. 
\subsection{Convergence of Nonlinear Subdivision Schemes in $L^p$ Spaces}
In the sequel,  we give a sufficient 
condition for the convergence of nonlinear subdivision 
schemes in $L^p(\R^{d})$. This result will be a generalization of the 
existing result in the 
linear context established in \cite{HanJia} and only uses the operator $S_1$. 
\newtheorem{ConvContrac}{Theorem}
\begin{ConvContrac}
\label{ConvContract1}
Let $S$ be a data dependent subdivision operator that reproduces the 
constants. If $\rho_{p,1}(S) < m^{\frac{1}{p}}$, then  
$Sv^j$ converges to a $L^p$ limit function.
\end{ConvContrac} 
\textsc{Proof:} 
Let us consider 
\begin{equation}
\label{fjdef}
v_j(x) := \sum\limits_{k \in \bZ} v^{j}_k 
\phi_{j,k}(x), 
\end{equation}
where 
$\phi(x) = \prod\limits_{i=1}^{d} \max(0,1-|x_i|)$ is the hat function. 
With this choice, one can
easily check that $\sum\limits_{k \in \bZ} 
\phi(x-k) = 1$. 
Let  $\rho_{p,1} (S) < \rho < m^{\frac{1}{p}}$,
it follows that 
$$
\|\Delta^1 v^j\|_{(\lp)^d} \lsim \rho^j \|\Delta^1 v^0\|_{(\lp)^d},
$$
since $q_1=d$.
We now show that the sequence $v_j$ is a Cauchy sequence in $L^p$:
\begin{eqnarray*}
v_{j+1}(x) - v_j(x) 
= 
\sum\limits_{k \in \bZ} 
v^{j+1}_k 
\phi_{j+1,k}(x) - 
\sum\limits_{p \in \bZ} 
v^{j}_p 
\phi_{j,p}(x)
\\
= 
\sum\limits_{k \in \bZ} 
\sum\limits_{p \in \bZ} 
(v^{j+1}_k-v^{j}_p)  
\phi_{j+1,k}(x)
\phi_{j,p}(x) 
\end{eqnarray*}
where we have used $\sum\limits_{k \in \bZ} \phi(\cdot-k)= 1 .$
Now, since the subdivision operator reproduces the constants: 
\begin{eqnarray*}
v_{j+1}(x) - v_j(x) 
=\sum\limits_{p \in \bZ}
\sum\limits_{k \in \bZ}
\sum\limits_{l \in \bZ}
a_{k -M l}(v^j) (v^{j}_l-v^{j}_p)  
\phi_{j+1,k}(x)
\phi_{j,p}(x).
\end{eqnarray*}
Note that 
\begin{eqnarray*}
\sum\limits_{l \in \bZ}
a_{k -M l}(v^j) (v^{j}_l-v^{j}_p)  = 
\sum\limits_{l \in \bZ} a_{k -M l}(v^j) v^{j}_ l - v^j_p =   \sum\limits_{l \in \bZ} (a_{k -M l}(v^j) -\delta_{p-l}) v^{j}_ l
.\end{eqnarray*} 

Since  $\sum\limits_{l \in \bZ} a_{k -M l}(v^j) -\delta_{p-l} = 0$,  
$
\left \{ \nabla_i \delta_{l-\beta}, l \in \left \{ F(k)  
\cup \{p\}\right \}, \beta  \in \bZ, i=1,\cdots,d \right \}
$ spans $(a_{k-Ml}-\delta_{p-l})_{ l \in \left \{ F(k)  
\cup \{p\}\right \}}$.
This enables us to write:
\begin{eqnarray*}
v_{j+1}(x) - v_j(x) 
= 
\sum\limits_{p \in \bZ}
\sum\limits_{k \in \bZ}
\sum\limits_{l \in V(k)\bigcup \{p\}}
\sum\limits_{i=1}^d d^i_{k,p,l} \nabla_i v^j_l
\phi_{j+1,k}(x)
\phi_{j,p}(x). 
\end{eqnarray*}
Since $|\sum\limits_{k \in \bZ} 
\phi_{j+1,k}(x)| = 1$  following the same argument 
as in Theorem 3.2 of \cite{HanJia}, we may write:
\begin{eqnarray}
\label{diffj}
\|v_{j+1} - v_j\|_{\Lp} 
&\lsim&
m^{-\frac{j}{p}} \max\limits_{1 \leq i \leq d} 
\| \nabla_i v^j\|_{\lp}\nonumber\\
&\lsim& m^{-\frac{j}{p}}\|\Delta^1 v^j\|_{\lp}\nonumber\\
&\sim& (\frac{\rho}{m^{\frac{1}{p}}})^j \|\Delta^1 v^0\|_{\lp}
\end{eqnarray}
which proves that $v_j$ converges in $L^p$, since $\rho < m^{\frac{1}{p}}$. 
Note that, for $p=\infty$, we obtain that the limit function is continuous.
\fin
\\
Furthermore, the above proof is valid for any function $\Phi_0$ satisfying the property of partition of unity when $p=\infty$. 
In general, we could show, following Theorem 3.4 of  \cite{HanJia}, 
that the limit function in $L^p$ is independent of the choice 
of a continuous and compactly supported $\Phi_0$. 
%%%%%%%%%%%%%%%%%%%%%%%%%%%%%%%%%%%%%%%%%%%%%%%%%%%%%%%%%%%%%%%%%%%%%%%%%%
%%%%%%%%%%%%%%%%%%%%%%%%%%%%%%%%%%%%%%%%%%%%%%%%%%%%%%%%%%%%%%%%%%%%%%%%%
\subsection{Uniform Convergence of the Subdivision Schemes to $C^s$ functions 
$(s < 1)$}
We are now ready to establish a sufficient condition for the $C^s$ smoothness 
of the limit function with $s<1$.
\newtheorem{ConvContracS}[ConvContrac]{Theorem}
\begin{ConvContracS}
\label{ConvContrac2}
Let $S(v)$ be a data dependent subdivision operator which reproduces the constants.
If the scheme for the differences satisfies 
$\rho_{p,1} (S)< m^{-s+\frac{1}{p}}$, for some $0< s < 1$ then 
$Sv^j$ is convergent in $L^p$ and the limit 
function  is $C^{s}$ .
\end{ConvContracS} 
\textsc{Proof:}
First, the convergence in $L^p$ is a consequence of 
$\rho_{p,1} (S)< m^{\frac{1}{p}}$. 
In order to prove that the limit function $v$ be in $C^s$, it suffices to 
evaluate $|v(x)-v(y)|$ for $\| x - y\|_{\infty}\leq 1$. Let $j$ be
such that $m^{-j-1} \leq \| x-y\|_{\infty} \leq m^{-j}$.
We then write :
\begin{eqnarray*}
|v(x)-v(y)|
&\leq&|v(x)-v_j(x)|+|v(y)-v_j( y)|+
|v_j(x)-v_j(y)| \\
&\leq& 2\|v-v_j\|_{\Li}+|v_j(x)-v_j(y)|
\end{eqnarray*}
Note that (\ref{diffj}) implies that 
$\|v-v_j\|_{\Li}\lsim \rho^j \|\Delta^1 v^0\|_{\li}$. 
Since $v_j$ is absolutely continuous, it is almost everywhere 
differentiable, so putting $y = x+M^{-j}h$,
with $h = (h_i)_{i=1,\cdots,d}$ satisfying $\|h\|_\infty \leq 1$ we get:
\begin{eqnarray*}
|v_j(x+ M^{-j}h)-v_j(x)|
&\leq&
|v_j(x+M^{-j}h)-v_j(x+M^{-j}(h-h_d e_d))| \\
&&+|v_j(x+ M^{-j}(h-h_d e_d))-
 v_j(x+M^{-j}(h-h_d e_d-
 h_{d-1} e_{d-1}))|\\
&&+\cdots+
|v_j(x+M^{-j}(h_1 e_1))-
 v_j(x)|
\end{eqnarray*}
Then, using a Taylor expansion we remark that, there exists $\theta \in
]-h_d,h_d[$ such that:
\begin{eqnarray*}
|v_j(x+M^{-j}h)-v_j(x+M^{-j}(h-h_d e_d))| &=&
\sum_{k \in \bZ} 
v_k^j h_d D_d\phi(M^{j}x+h-k+\theta_d e_d) 
\end{eqnarray*}
If we denote $\Psi_d(x) = \Phi(x_1)\cdots \Phi(x_{d-1})\Psi(x_d)$, where 
$\Psi$ is the characteristic function of $[0,1]$ and $\Phi(x_i) =
\max(0,1-|x_i|)$, we may write:
\begin{eqnarray*}
|v_j(x+M^{-j}h)-v_j(x+M^{-j}(h-h_d e_d))| &\sim&
\sum_{k \in \bZ} \nabla_d v_k^j h_d \Psi_d(y) 
\end{eqnarray*}
where $y_i=(M^{-j}x+h-k+\theta_d e_d)_i$ if $i <d$ and $y_d = 2(M^{-j}x+h+
\theta_d e_d)_d-k_d$ (we have used the fact that the differential of the hat
function $\Phi$ is the Haar wavelet). 
Iterating the procedure for other differences in the sum, we get:
\begin{eqnarray*}
|v_j(x+ M^{-j}h)-v_j(x)|
&\lsim&  \sum_{i=1}^d \|\nabla_i v^{j}\|_{\li}  
 \lsim  \|\Delta^1 v^j\|_{\li}.
\end{eqnarray*}
Combining these results we may finally write:
\begin{eqnarray*}
|v(x) -v(y)|&\leq& |v(x)-v_j(x)|
                          +|v(y)-v_j(y)|
             	  +|v_j(x)-v_j(y)|\\
&\leq&2\|v-v_j\|_{\Li} + |v_j(x) -v_j(y)|\\
&\lsim&\rho^j\|\Delta^1 v^0 \|_{\li}+ \|\Delta^1 v^j\|_{\li}\\
&\lsim&\rho^j \|\Delta^1 v^0 \|_{\li}\lsim \|x-y\|^s_{\infty}
\end{eqnarray*}   
with $s < -\log(\rho_{\infty,1} )/\log m$.
\fin 
\section{Examples of Bidimensional Subdivision Schemes }

In the  first part of this section, we construct an interpolatory subdivision scheme having the  dilation matrix  the hexagonal matrix which is: 
\begin{eqnarray*}
{  M} &=& \left (
\begin{array}{c c}
2&1\\
0&-2
\end{array}
\right ),
\end{eqnarray*}
For the hexagonal dilation matrix, the coset vectors are $\varepsilon_0=(0,0)^T
,\varepsilon_1=(1,0)^T,\varepsilon_2=(1,-1)^T,\varepsilon_3=(2,-1)^T$.
The coset vector $\varepsilon_i$,$i=0,\cdots,3$ of $M$ defines a partition of $\mathbb{Z}^2$ as 
follows:
$$
\mathbb{Z}^2 = \bigcup\limits_{i=0}^{3} \left 
\{ Mk+\varepsilon_i, k \in \mathbb{Z}^2 \right \}.
$$

The discrete data at the level $j$,  $v^j$ is defined on the grid $\Gamma^{j}= M^{-j}\mathbb{Z}^2$, the value $v_k^j$ is then associated to the location 
$M^{-j}k$.
We now define  our  bi-dimensional 
interpolatory subdivision scheme based on 
the data dependent  subdivision operator which  acts from the coarse 
grid $\Gamma^{j-1}$ to the fine grid grid $\Gamma^{j}$.
To this end,  we will compute $v^j$ at the different coset points on the fine grid $\Gamma^{j}$ using  the existing values $v^{j-1}$ of the coarse grid $\Gamma^{j-1}$, 
as follows: for the first coset vector 
$\varepsilon_0=(0,0)^T$ we simply put 
$v^j_{Mk +\varepsilon_0 } = v^{j-1}_k$, 
for the coset vectors ${\varepsilon_i}$,  $i=1,\cdots,3$.
the value $v^j_{Mk+\varepsilon_i}$,  $i=1,\cdots,3$ 
is defined by affine interpolation of the values on the coarse grid.
 To do so, 
we define four different stencils on $\Gamma^{j-1}$ as follows:
\begin{eqnarray*}
V_k^{j,1} &=& \{M^{-j+1}  k, M^{-j+1}(k+e_1),M^{-j+1}(k+e_2)\},\\
V_k^{j,2} &=& \{M^{-j+1} k, M^{-j+1}(k+e_2),M^{-j+1}(k+e_1+e_2)\},\\
W_k^{j,1} &=& \{M^{-j+1} (k+e_1),M^{-j+1}(k+e_2),M^{-j+1}(k+e_1+e_2)\},\\
W_k^{2} &=& \{M^{-j+1} k,M^{-j+1}(k+e_1),M^{-j+1}(k+e_1+e_2)\}.
\end{eqnarray*}
We determine  to which stencils each point of $\Gamma^j$ belongs to, and we then define the prediction as its barycentric coordinates. Since we use an affine interpolant we have: 
\begin{eqnarray}
\label{lb1}
v^{j}_{Mk} = v^{j-1}_k \textrm{ and }
v^{j}_{Mk + \varepsilon_1} = \frac{1}{2} v^{j-1}_{k} +
\frac{1}{2} v^{j-1}_{k+e_1}. 
\end{eqnarray}
To compute the rules for the coset point $\varepsilon_2$,
$V_k^{1}$ or $V_k^{2}$ can be used  leading respectively to:
\begin{eqnarray}
\label{lb2}
v^{j,1}_{Mk+\varepsilon_2} &=& \frac{1}{4} v^{j-1}_{k+e_1} +
\frac{1}{2} v^{j-1}_{k+e_2}
                +\frac{1}{4} v^{j-1}_{k} \nonumber\\
v^{j,2}_{Mk+\varepsilon_2} &=& \frac{1}{2} v^{j-1}_k +
                      \frac{1}{4} v^{j-1}_{k+e_2} +
                      \frac{1}{4} v^{j-1}_{k+e_1+e_2}.
\end{eqnarray}
When one considers the rules for the coset point $\varepsilon_3$,
$W_k^{1}$ or $W_k^{2}$ can be used leading respectively to:
\begin{eqnarray}
\label{lb3}
v^{j,1}_{Mk+\varepsilon_3} &=& \frac{1}{4} v^{j-1}_{k+e_2} +
                      \frac{1}{4} v^{j-1}_{k+e_1+e_2}
                     +\frac{1}{2} v^{j-1}_{k+e_1} \nonumber \\
v^{j,2}_{Mk+\varepsilon_3} &=& \frac{1}{4} v^{j-1}_k +
                      \frac{1}{4} v^{j-1}_{k+e_1} +
                      \frac{1}{2} v^{j-1}_{k+e_1+e_2}.
\end{eqnarray}
This nonlinear scheme is converging in $L^\infty$ since we have the following result: 
\newtheorem{conver}[DataDependent]{Proposition}
\begin{conver}
The prediction defined by (\ref{lb1}), (\ref{lb2}), (\ref{lb3}) satisfies:
$$
\|\Delta^1 v^j_{M.+\varepsilon_i}\|_{\li} \leq \frac{3}{4} 
\|\Delta^1 v^{j-1}\|_{\li}
$$
\end{conver}
We do not detail the proof here but the result is obtained by computing the differences in the canonical directions 
at each coset points. 
If we then use Theorem \ref{ConvContrac2} we can find the regularity of the 
corresponding limit function 
is $C^s$ with $s< -\frac{\log (3/4)}{\log(4)} \approx 0.207$.

In the second pat of the section, we build an 
example of bidimensional subdivision scheme based on the same 
philosophy but this time using as dilation matrix the quincunx matrix 
defined by:
\begin{eqnarray*}
M&=&\left (
\begin{array}{cc}-1 & 1\\
1&1
\end{array}
\right ),
\end{eqnarray*}
whose coset vectors are $\varepsilon_0=(0,0)^T$ and
$\varepsilon_1=(0,1)^T$.
Note that $a_{0,0}=1$ and since the nonlinear subdivision operator
reproduces the constants we have
$\sum\limits_{i} a_{Mi+\varepsilon}=1$ for all coset vectors
$\varepsilon$. To build the subdivision operator, we consider the subdivision rules based on interpolation 
by of first degree polynomials on the grid $\Gamma^{j-1}$. 
$v^j_{Mk+\epsilon_1}$  corresponds to a point 
inside the cell delimited by 
$M^{-j+1} \{ k,k+e_1,k+e_2,k+e_1+e_2\}$. 
There are four potential stencils,  leading in this case only to  two subdivision rules:
\begin{eqnarray}
\label{quin1}
\hat v^{j,1}_{Mk+\epsilon_1}&=\frac 1 2 (v_{k}^{j-1}+v_{k+e_1+e_2}^{j-1})
\end{eqnarray}
\begin{eqnarray}
\label{quin2}
\hat v^{j,2}_{Mk+\epsilon_1}&=\frac 1 2 (v_{k+e_1}^{j-1}+v_{k+e_2}^{j-1})
\end{eqnarray}
Note also, that since the scheme is interpolatory we have the
relation: $v^j_{Mk} = v^{j-1}_k$. Let us 
now prove a contraction property for the above scheme.

\newtheorem{convquin}[DataDependent]{Proposition}

\begin{convquin}
\label{convquin1}
The nonlinear subdivision scheme defined by (\ref{quin1}) and (\ref{quin2}) satisfies the following property:
\begin{enumerate}
\item when $k=Mk'$:
\begin{eqnarray*}
\|v^{j,1}_{M.+\varepsilon_1}-v^j_{M.}\|_{\li}&\leq& 
\frac{1}{2}\|\Delta^1 v_.^{j-2}\|_{\li}\\
\|v^{j,2}_{M.+\varepsilon_1}-v^j_{M.}\|_{\li}&\leq& 
\|\Delta^1 v_.^{j-2}\|_{\li}
\end{eqnarray*}
\item when $k=Mk'+\varepsilon_1$, we can show that:
\begin{eqnarray*}
\|v^{j,2}_{M.+\varepsilon_1}-v^j_{M.}\|_{\li}&\leq& \frac{1}{2} 
\|\Delta^1 v_.^{j-2}\|_{\li}\\
\|v^{j,1}_{M.+\varepsilon_1}-v^j_{M.}\|_{\li}&\leq& 
\|\Delta^1 v_.^{j-2}\|_{\li}
\end{eqnarray*}
\end{enumerate}
\end{convquin}
The proof of this theorem is obtained computing all the potential differences.
This theorem shows that the nonlinear subdivision scheme converges in 
$L^\infty$ since $\rho_{1,\infty}(S) < 1$. 
\section{Convergence in Sobolev Spaces}
In this section, we extend the result established in \cite{Jia2} 
on the convergence of linear subdivision scheme to 
our nonlinear setting. We will first recall the notion of 
convergence in Sobolev spaces in the linear case. 
Following \cite{Jia} Theorem 4.2, when 
$\Phi_0(x) = \sum\limits_{k\in \bZ} a_k \Phi_0(Mx-k)$ is $L^p$-stable, the 
so-called "moment condition of order $k+1$ for $a$" is equivalent to the 
polynomial reproduction property of polynomial of total degree $k$ 
for the subdivision scheme associated to $a$. In what follows, we will say 
that $\Phi_0$ reproduces polynomial of total degree $k$. When the subdivision
associated to $a$ exactly reproduces polynomials, we will say that $\Phi_0$
exactly reproduces polynomials. 
We then have the following definition for the convergence of subdivision 
schemes in Sobolev spaces in the linear case \cite{Jia2}: 
\newtheorem{LinearSobo}[LocalPred]{Definition}
\begin{LinearSobo}
We say that $v^j=Sv^{j-1}$ converges in the Sobolev space $W_N^k(\bR)$ if there exists 
a function $v$ in $W_N^k(\bR)$ satisfying:
$$
\lim\limits_{j\rightarrow +\infty} \|v_j-v\|_{W_N^k(\bR)} = 0
$$
where $v$ is in $W_N^k(\bR)$, and $v_j=\sum\limits_{k \in \bZ} v_k^j \Phi_0(M^jx-k)$ 
for any $\Phi_0$ reproducing polynomials of total degree $k$.
\end{LinearSobo}
We are going to see that in the nonlinear case, to ensure the convergence 
we are obliged to make a restriction on the choice of $\Phi_0$.
We will first give some results when the matrix $M$ is an isotropic dilation 
matrix, we will also emphasize a particular class of isotropic matrices, 
very useful in image processing.
\subsection{Definitions and Preliminary Results}
\newtheorem{isotropic}[LocalPred]{Definition}
\begin{isotropic}
\label{isotropic1}
We say that a matrix $M$ is isotropic if it is similar to the diagonal matrix $\textrm{diag}(\sigma_1, \ldots, \sigma_d)$, i.e. there exists an invertible matrix $\Lambda$ such that
 \begin{equation*}
 M=\Lambda^{-1} \textrm{diag}(\sigma_1, \ldots, \sigma_d) \Lambda,
 \end{equation*}
 with $|\sigma_1| = \ldots = |\sigma_d|$ being the eigenvalues of matrix $M$.
\end{isotropic}
Evidently, for an isotropic matrix holds $|\sigma_1| = \ldots = |\sigma_d|=\sigma=m^{\frac{1}{d}}$. 
Moreover, for any given norm in $\mathbb{R}^{d}$, any integer 
$n$ and  any $v\in \mathbb{R}^{d}$ we have 
$$
\sigma^{n}\|u\| \lsim \|M^n u\| \lsim \sigma^{n} \|u\|.
$$
A particular class of isotropic matrices is when there exists a set 
$\tilde e_1, \tilde e_2,\cdots,\tilde e_q$ such that:
\begin{eqnarray}
\label{caspart}
M \tilde e_i = \lambda_i \tilde e_{\gamma(i)}
\end{eqnarray}
where $\gamma$ is a permutation of $\{1,\cdots,q\}$. Such matrices are 
particular cases of isotropic matrices 
since $M^q=\lambda I$ where $I$ is the identity matrix and where $\lambda = \prod\limits_{i=1}^d \lambda_i$. 
For instance, when $d=2$, the quincunx (resp. hexagonal) matrix satisfies $M^2=2I$ (resp. $M^2=4I$). 

We establish the following property on joint spectral radii that will be useful when dealing with the convergence in Sobolev spaces.
\newtheorem{spectral1}[DataDependent]{Proposition}
\begin{spectral1}
\label{spectral}
Assume that $S$ reproduces polynomials up to total 
degree $N$. Then,
$$
\rho_{p,n+1}(S)\geq
\frac {1} {\|M\|_{\infty}} \rho_{p,n}(S),
$$
for all  $n=0,\ldots,N$.
\end{spectral1}
\noindent \underline{Remark:}
If $M$ is an isotropic matrix and $S$ reproduces polynomials
up to total degree $N$, then
\begin{equation*}
\rho_{p,n+1}(S)\geq
\sigma^{-1} \rho_{p,n}(S),
\end{equation*}
for all  $n=0,\ldots,N$.\\
\textsc{Proof:} It is enough to prove
\begin{equation*}
\rho_{p,1}( S)\geq
\frac {1} {\|M\|_{\infty}} \rho_p(S).
\end{equation*}
According to the definition of spectral radius there exists
$\rho>\rho_{p,1}(S)$ such that for any $u^0$
\begin{equation*}
\|S_{1}(u^{j-1}) \ldots S_{1}(u^0){\nabla} u\|_{\lp}
\lsim \rho^j \|{\nabla} u\|_{(\lp}.
\end{equation*}
Using the notation $\omega^j := S(u^{j-1})\cdot \ldots \cdot S(u^0) u$ 
we obtain
\begin{equation*}
\|{\nabla} \omega^j\|_{\lp} \lsim \rho^j \|{\nabla}
u\|_{\lp}.
\end{equation*}
Since
\begin{equation*}
\omega^j_l = \sum_n A^j_{l,n} u_n,
\end{equation*}
where
\begin{equation*}
A^j_{l,n} = \sum_{l_1, \ldots, l_{j-1}}a_{l-Ml_{j-1}}(u^{j-1})
a_{l_{j-1}-Ml_{j-2}}\cdot \ldots \cdot a_{l_1-Mn}(u^0).
\end{equation*}
We can write down the $\ell^p$-norm as follows:
\begin{equation*}
\|\omega^j\|^p_{\lp} =
\sum_{k\in \mathbb{Z}^d}\sum_{i=1}^{m^j}|\omega^j_{M^jk+\varepsilon^j_i}|^p,
\end{equation*}
where $\{\varepsilon^j_i\}_{i=1}^{m^j}$ are the representatives of cosets of $M^j$.
First note that:
\begin{equation*}
\|k-n\|_{\infty} \leq \|k-n+M^{-j}\varepsilon^j_i\|_{\infty}+
\|M^{-j}\varepsilon^j_i\|_{\infty}.
\end{equation*}
Note that $M^{-j} \varepsilon_i^j$ belongs to the unit square so that
$\|M^{-j}\varepsilon^j_i\|_{\infty} \leq K_1$.
When $A^j_{M^jk+\varepsilon^j_i,
n} \neq 0$, one can prove that there exists $K_2>0$ such that
\begin{equation*}
\|k-n+M^{-j}\varepsilon^j_i\|_{\infty} \leq K_2,
\end{equation*}
the proof being similar to that of Lemma 2 in \cite{Goodman}.
From these inequalities it follows that if
$A^j_{M^jk+\varepsilon^j_i, n} \neq 0$ there exists $K_3>0$ such that
\begin{equation*}
\|k-n\|_{\infty} \leq K_3,
\end{equation*}
that is, for a fixed $k$, the values of $\omega^j_l$ for
$l \in \{M^jk+\varepsilon^j_i\}^{m^j}_{i=1}$ depend only on
$u_n$ with $n : \{\|k-n\|_{\infty} \leq K_3\}$.

\noindent Let us now fix $k$ and define $\tilde{u}$ such that
$$
\tilde{u}_l=\begin{cases}
u_l,&\text{if $\|k-l\|_{\infty} \leq K_3$;}\\
0,&\text{otherwise.}
\end{cases}
$$
Let $\tilde{\omega}^j := S(u^{j-1})\cdot \ldots \cdot S(u^0) \tilde{u}$, then
$$
\tilde{\omega}^j_l=\begin{cases}
\omega^j_l,&\text{ if $l \in \{M^jk+\varepsilon^j_i\}^{m^j}_{i=1}$;}\\
0,&\text{ if $\|k-M^{-j}l\|_{\infty} \geq K_4$,}
\end{cases}
$$
since if $A^j_{l, n}\neq 0$, then
\begin{equation*}
\|k-M^{-j}l\|_{\infty} \leq \|k-n\|_{\infty} + \|n-M^{-j}l\|_{\infty} \leq K_3 +
K_2  :=  K_4.
\end{equation*}
Moreover, from $\|k-M^{-j} l\|_{\infty} \leq K_4$, it 
follows that $\|M^{j}k-l\|_{\infty} \leq K_4 \|M^{j}\|_{\infty}$.
 Taking all this into account, we get
\begin{eqnarray*}
\sum\limits_{k \in \Z^d}
\sum_{l \in \{M^jk+\varepsilon^j_i\}}
|\omega^j_l|^p &=&
\sum\limits_{k \in \Z^d} \sum_{l \in \{M^jk+\varepsilon^j_i\}}|
\tilde{\omega}^j_l|^p \leq
\sum\limits_{k \in \Z^d}
\sum_{\|M^{j}k-l\|_{\infty} \leq K_4 \|M^{j}\|_{\infty}}
|\tilde{\omega}^j_l|^p \\
&\lsim& \|M\|_{\infty}^j \|\Delta^1 \tilde{\omega}^j_l\|_{\lp}
\lsim (\|M\|_{\infty} \rho)^j \|\Delta^1 \tilde{u}\|_{\lp}.
\end{eqnarray*}
That is, $\|\omega^j\|_{\lp} \lsim (\|M\|_{\infty} \rho)^j
\|u\|_{\lp}$, consequently $\rho_p(S) \lsim \|M\|_{\infty} \rho$.
Now, if $\rho \rightarrow \rho_{p,1}(S)$ we get
$\rho_p(S) \leq \|M\|_{\infty} \rho_{p,1}( S)$.
\fin
\subsection{Convergence in Sobolev Spaces When $M$ is Isotropic} 
First, Let us recall that the Sobolev norm on $W_N^p(\bR)$ is defined by:
\begin{eqnarray}
\label{defSobo}
\|f\|_{W_N^p(\bR)} = \|f\|_{\Lp} + \sum\limits_{|\mu| \leq N} \|D^\mu f\|_{\Lp}.
\end{eqnarray}
If one considers a set $x_{\bf 1},\cdots,x_{\bf n}$ such that $
[x_{\bf 1},\cdots,x_{\bf n}]\mathbb{Z}^n=\mathbb{Z}^d$, an equivalent norm 
is given by:
\begin{eqnarray}
\label{defSobo1}
\|f\|_{W_N^p(\bR)} = \|f\|_{\Lp} + \sum\limits_{|\mu| \leq N} \|\tilde D^\mu f\|_{\Lp}.
\end{eqnarray}
where 
$\tilde D^{\mu} = D^{\mu_1}_{x_{\bf 1}} \cdots D^{\mu_n}_{x_{\bf n}}$.
\\
We then enounce a convergence theorem for general isotropic matrix $M$:
\newtheorem{segaln}[ConvContrac]{Theorem}
\begin{segaln}
\label{segaln}
Let $S$ be a data dependent nonlinear subdivision scheme which exactly 
reproduces polynomials up to total degree $N-1$, then the subdivision scheme 
$Sv^j$ converges in $W_N^p(\R^{d})$, provided 
$\Phi_0$ is compactly supported and exactly reproduces polynomials up to total degree $N-1$ and    
\begin{eqnarray}
\label{rspect1}
\rho_{p,N}(S) < m^{\frac{1}{p}-\frac{s}{d}} \textrm{ for some }s > N. 
\end{eqnarray}
\end{segaln}
\textsc{ Proof:}
Note that because of Proposition \ref{spectral}, the hypotheses of Theorem 
\ref{segaln} imply 
that  
$
\rho_{p,k}(S)  < m^\frac{1}{d} \rho_{p,k+1}(S) < m^{\frac{1}{p}-\frac{s-1}{d}},
$
which means that (\ref{rspect1}) is also true for $k < N$.
Let us now show that $v_j$ is a Cauchy sequence in $L^p$. To do so, let us define  
$$
q_j(x)=\sum_{l=1}^d \lambda_{j,l}x_l,
$$
where $\Lambda = (\lambda_{j,l})$ is defined in (\ref{isotropic1}).
For a multi-index $\mu=(\mu_1, \ldots, \mu_d) \in \mathbb Z^d$ let
$$
q_\mu (x) = q_1^{\mu_1}(x) \ldots q_d^{\mu_d}(x).
$$
Since $\Lambda$ is invertible, the set $\{q_{\mu}: |\mu| = N\}$ 
forms a basis of the space of all polynomials of exact degree $N$, 
which proves that
$$
\|D^\mu (v_{j+1}-v_j)\|_{\Lp} \sim
\|q_{\mu}(D)(v_{j+1}-v_j)\|_{\Lp}
$$
Now, we use the fact that, since $M$ is isotropic,  
$
q_\mu(D) (f(M^jx))=\sigma^{j\mu}(q_\mu(D) f)(M^jx)
$
where $\tilde \sigma^\mu = \prod\limits_{i=1}^d \sigma_i^{\mu_i}$ 
(\cite{Jia1}). We can thus write:
\begin{eqnarray*}
q_\mu(D)(v_{j+1}-v_j)=
q_\mu(D) \left(
\sum_{l \in \bZ} v^{j+1}_l \Phi_0 (M^{j+1}x-l) - 
\sum_{l \in \bZ}  
v^j_l  \Phi_0 (M^jx-l)\right ). 
\end{eqnarray*}
We use now the scaling equation of $\Phi_0$ to get 
\begin{eqnarray*}
q_\mu(D)(v_{j+1}-v_j)
&=&q_\mu(D) 
\left(
\sum_{l \in \bZ}v^{j+1}_l  \Phi_0 (M^{j+1}x-l) - 
\sum_{l \in \bZ}  
\sum_{r \in \bZ} v^j_r g_{l-Mr} \Phi_0 (M^{j+1} x- l) \right )\\
&=& \sum_{l \in \bZ} (v^{j+1}_l - \sum_{r \in \bZ} v^j_r g_{k-Mr} ) 
q_\mu (D) \left ( \Phi_0 (M^{j+1}x-l) \right )\\
&=& \sum_{l \in \bZ}  \sum_{r \in \bZ} (a_{l-Mr}(v^j)-g_{l-Mr}) v^j_r 
\tilde \sigma^{\mu(j+1)}  (q_{\mu}(D) \Phi_0) (M^{j+1}x-l). 
\end{eqnarray*}
Since $S$ and $\Phi_0$ exactly reproduce polynomials up to total degree 
$N-1$, we have for $|\mu| \leq N-1$:
$$
\sum_{r \in \bZ} (a_{l-Mr}(v^j)-g_{l-Mr}) r^{\mu} = 0.
$$
Remark that $g_{l-Mr}= 0$ for $\|l-Mr\| > \tilde K$ since $\Phi_0$ is compactly
supported.
Since 
$
\left \{
\nabla^{\nu} \delta_{l-\beta}, |\nu| = N, r 
\in F(l)=\left \{ \|l-Mr \| \leq \max(K,\tilde K) \right \}, \beta \in \bZ 
\right \}
$ spans $(a_{l-Mr}(v^j)-g_{l-Mr})_{r \in F(l)}$, 
we deduce:
\begin{eqnarray*}
q_\mu(D)(v_{j+1}-v_j) &=& 
\sum_{l \in \bZ} \sum_{r \in F(l)} \sum_{|\nu|=N}c^{\nu}_{r}(v^{j})
\nabla^{\nu} v^{j}_r \tilde \sigma^{\mu (j+1)}  (q_{\mu}(D) \Phi_0) (M^{j+1}x-l), 
\end{eqnarray*}
Consequently,
\begin{eqnarray*}
\|q_{\mu}(D)(v_{j+1}-v_j)\|_{\Lp} 
&\lsim& \sigma^{(j+1)N} m^{-(j+1)/p} (\rho_{p,N}(S))^j 
\|\Delta^{N} v^{0}\|_{(\lp)^{q_{N}}}
\end{eqnarray*}
Since $\rho_{p,N}(S) < m^{1/p-s/d}$, with $s > N$ we obtain
\begin{eqnarray*}
\|q_{\mu}(D)(v_{j+1}-v_j)\|_{\Lp} &\lsim&  \sigma^{j(N-s)} 
\|\Delta^{N} v^0\|_{(\lp)^{q_{N}}}.
\end{eqnarray*}
From this we deduce that $\|q_{\mu}(D)(v_{j+1}-v_j)\|_{\Lp}$ tends to $0$ with $j$. Making $\mu$ vary, we deduce the 
convergence in $W_{N}^p(\bR)$ \fin
 
We now show that when the matrix $M$ satisfies (\ref{caspart}) and when $\Phi_0$ is a box spline 
satisfying certain properties, the limit function is in $W_N^p(\bR)$. Before that, we need to recall 
the definition of box splines and some properties that we will use.   
Let us define a set of $n$ vectors, not necessarily distinct:
$$
X_n = \{ x_{\bf 1}, \cdots,x_{\bf n} \} 
\subset \bZ \setminus \{ 0 \}.
$$
We assume that $d$ vectors of $ X_n$ are linearly independent.
Let us rearrange the family $X_n$ such that 
$X_d = \{ x_{\bf 1},\cdots,x_{\bf d} \}$ are linearly
independent. We denote by
$[x_{\bf 1},\cdots,x_{\bf d}][0,1[^{d}$ the collection of linear combinations 
$\sum\limits_{i=1}^{d} \lambda_i x_{\bf i}$ with $\lambda_i \in [0,1[$.
Then, we define multivariate box splines as follows \cite{Ch}\cite{Pr}:
\begin{eqnarray}\label{defSpli}
\beta_0(x,X_d ) &=&\left \{ \begin{array}{c}
              \frac{1}{|\det(x_{\bf 1},\cdots,x_{\bf d})|} \textrm{ if } 
	      x \in [x_{\bf 1},\cdots,x_{\bf d}][0,1[^{d}\\
	      0 \textrm{ otherwise}
	      \end{array}
	      \right . \nonumber \\ 
\beta_0(x,X_k) &=& \int_0^1 \beta_0 (x-tx_{\bf k},X_{k-1}) dt, 
\quad n \geq k > d.
\end{eqnarray}
One can check by induction that the support of $\beta_0(x,X_n)$ is 
$[x_{\bf 1},x_{\bf 2},\cdots,x_{\bf n}][0,1]^n$. The regularity of box
splines is then given by the following theorem \cite{Pr}:
\newtheorem{Differentiability}[DataDependent]{Proposition}
\begin{Differentiability}
\label{Different}
$\beta_0(x,X_n)$ is $r$ times continuously differentiable if all
subsets of $X_n$ obtained by deleting $r+1$ vectors spans
$\bR$.
\end{Differentiability}
We recall a property on the directional derivatives 
of box splines, which we use in the convergence theorem that follows:  
\newtheorem{Differentiation}[DataDependent]{Proposition}
\begin{Differentiation}
\label{Different1}
Assume that $X_n \setminus x_{\bf r}$ spans $\bR$, and 
consider the following box spline function $s(x) = \sum\limits_{k \in \bZ}
c_k \beta_0(x-k,X_n)$  
then the directional derivative of $s$ in the direction $x_{\bf r}$ 
reads:
$$
D_{x_{\bf r}} s(x) = \sum\limits_{k \in \bZ} 
\nabla_{x_{\bf r}} c_k 
\beta_0(x-k,X_n \setminus x_{\bf r}).
$$
\end{Differentiation}  
We will also need the property of polynomial reproduction which is \cite{Pr}:     
\newtheorem{Polrep}[DataDependent]{Proposition}
\begin{Polrep}
\label{Polrep1}
If $\beta_0(x,X_n)$ is r times continuously differentiable then, for any polynomial $c(x)$ of total degree $d \leq r+1$,
\begin{eqnarray}
\label{defv}
p(x) = \sum_{i \in \bZ} c(i) \beta_0(x-i,X_n)
\end{eqnarray}
is a polynomial with total degree $d$, with the same leading coefficients
(i.e. the coefficients corresponding to degree $d$). Conversely, for any 
polynomial $p$, it satisfies (\ref{defv}) with 
$c$ being a polynomial having the same leading coefficients as $p$.
\end{Polrep}
\newtheorem{segaln2}[ConvContrac]{Theorem}
\begin{segaln2}
\label{segaln1}
Let $S$ be a data dependent nonlinear subdivision scheme which reproduces 
polynomials up to total degree $N-1$ and assume that $M$ satisfies 
relation (\ref{caspart}), then the subdivision scheme $Sv^j$ converges 
in $W_N^p(\R^{d})$, if when $N \geq 2$, $\Phi_0$ is a 
$C^{N-2}$ box spline generated by $x_{\bf 1},\cdots,x_{\bf n}$ 
satisfying $\Phi_0(x) =\sum\limits_{k} g_k \Phi_0(Mx-k)$ and if
$N =1$  $\Phi_0(x) =\sum\limits_{k} g_k \Phi_0(Mx-k)$ and 
$\sum\limits_{k \in \bZ} \Phi_0(x-k)=1$ and if 
\begin{eqnarray}
\label{rspect}
\rho_{p,N}(S)  <  m^{\frac{1}{p}-\frac{s}{d}} \textrm{ for some } s > N. 
\end{eqnarray}
\end{segaln2}
\textsc{Proof:} We here prove the case $N \geq 2$, the case $N=1$ can be proved
similarly. First note that since $\Phi_0(x)$ is a $C^{N-2}$ box spline, 
we can write for any polynomial $p$ of total degree $N-1$ at most:
\begin{eqnarray*}
p(M^{-1}x) &=& \sum\limits_{i \in \bZ} \tilde p(i) \Phi_0(M^{-1}x-i,X_n)\\
&=&\sum\limits_{q\in \bZ} \sum\limits_{i \in \bZ} g_{q-Mi} \tilde p(i) 
\Phi_0(x-q,X_n)
\end{eqnarray*}
Using Proposition \ref{Polrep1} we get $p$ and $\tilde p$ have the same leading coefficients, and that 
$\sum\limits_{i \in \bZ} g_{q-Mi} \tilde p(i)$ is a polynomial evaluated in 
$M^{-1}i$ having the same leading 
coefficients as $p$. That is to say the subdivision scheme 
$(Sv^j)_q = \sum\limits_{i \in \bZ} g_{q-Mi} v^j_i$ reproduces 
polynomials up to degree $N-1$.
 
As already noticed, the joint spectral radius of difference operator is independent of the choice of the 
directions $x_{\bf 1},\cdots,x_{\bf n}$ that spans $\bZ$. 
Furthermore, it is shown  in \cite{Me}, that the existence of a scaling 
equation for $\Phi_0$ implies that the vectors $x_{\bf i}$, ${\bf i} = 
{\bf 1},\cdots,{\bf n}$ satisfy a relation of type (\ref{caspart}).   
We consider such a set $\{ x_{\bf i} \}_{{\bf i}={\bf 1},\cdots,{\bf n}}$ 
and then define $\Phi_0(x)=\beta_0(x,Y_N)$ the box spline associated 
to the set 
$$
Y_N:=\left \{ \overbrace{x_{\bf 1},\cdots,x_{\bf 1}}^{N},\cdots,
\overbrace{x_{\bf n},\cdots,x_{\bf n}}^{N} \right \}.
$$
which is $C^{N-2}$ by definition. We then define  the differential operator 
$ \tilde D_{M^{-j}}^\mu :=\tilde D^{\mu_1}_{M^{-j}
x_{\bf 1}}\cdots \tilde D^{\mu_n}_{M^{-j} x_{\bf n}}$. We will use the
characterization (\ref{defSobo1}) of Sobolev spaces therefore $\mu =
(\mu_i)_{i=1,\cdots,n}$. For any $|\mu| \leq N$ we may write:
\begin{eqnarray*}
\tilde D^{\mu}_{M^{-j-1}} (v_{j+1}(x) - v_j(x)) 
&=& \sum\limits_{k \in \bZ} 
v^{j+1}_k 
(\tilde D^{\mu} \beta_0)(M^{j+1}x-k,Y_N)\\ 
&&- 
\sum\limits_{p \in \bZ} \sum_{i \in \bZ}
v^{j}_i  g_{p-Mi} (\tilde D^{\mu} \beta_0)(M^{j+1}x - p,Y_N),\\
\end{eqnarray*} 
using the scaling property satisfied by $\beta_0$.
Then, we get:
\begin{eqnarray*}
\tilde D^{\mu}_{M^{-j-1}} (v_{j+1}(x) - v_j(x)) 
&=& \sum\limits_{k \in \bZ}  \sum_{i \in \bZ}
(a_{k-Mi}(v^j)-g_{k-Mi}) v^j_i (\tilde D^{\mu} \beta_0)(M^{j+1}x-k,Y_N) \\
&=&\sum\limits_{k \in \bZ} \tilde \nabla^{\mu} ( \sum_{i \in \bZ}
(a_{k-Mi}(v^j)-g_{k-Mi}) v^j_i) \beta_0(M^{j+1}x-k,Y_N^\mu) 
\end{eqnarray*} 
where $Y_N^\mu$ is obtained by removing $\mu_i$ vector $x_{\bf i}$, 
$i=1,\cdots,d$ to $Y_N$ and $\tilde \nabla^{\mu} = 
\left (\nabla^{\mu_i}_{x_{\bf i}} \right )_{i=1,\cdots,n}$.
As both $a_{k-M.}(v^j)$ and $g_{k-M.}$ reproduce polynomials up to total 
degree $N-1$, there exist a finite sequence $c_{k,p}$ such that: 
\begin{eqnarray*}
\tilde \nabla^{\mu} ( \sum_{i \in \bZ} (a_{k-Mi}(v^j)-g_{k-Mi}) v^j_i)
&=&\sum\limits_{p \in V(k)\bigcup \tilde V (k)} 
\sum\limits_{|\nu|=|\mu|} c_{k,p}(\nu)  \tilde \nabla^\nu v^{j}_p,
\end{eqnarray*}
where $\tilde V (k) = \{i,\|k-Mi\| \leq \tilde K\}$, where $g_{k-Mi} = 0$ 
if $\|k-Mi\| > \tilde K$.
We finally deduce that:
\begin{eqnarray*}
\tilde D_{M^{-j-1}}^{\mu} (v_{j+1}(x) -v_j(x))=
\sum\limits_{k \in \Z^{d}} \sum\limits_{p \in V(k)\bigcup \tilde V(k)}
\sum\limits_{|\nu|=|\mu|} 
c_{k,p}( \nu ) \tilde \nabla^{\nu} v_p^j \beta_0(M^{j+1}x-k,Y_N^\mu).
\end{eqnarray*}
From this, we conclude that:
\begin{eqnarray*}
\| \tilde D_{M^{-j-1}}^{\mu} (v_{j+1}(x) -v_j(x)) \|_{\Lp} \lsim 
 \rho_{p,|\mu|}(S)^j m^{-\frac{j+1}{p}} \|\tilde \Delta^{|\mu|} v^{j}_0\|_
 {(\lp)^{\tilde q_{|\mu|}}}. 
\end{eqnarray*}
Now, consider a sufficiently differentiable function 
$f$ and remark that $D_{M^{-j-1}x_{\bf 1}} f(x) = (D f)(x).M^{-j-1}x_{\bf 1}$, 
where $Df$ is the differential of the function $f$. 
We also note that $M^q = \lambda I$ which implies that 
$\lambda = \sigma^q$ and we then put $j+1 = q \times\lfloor \frac{j+1}{q} 
\rfloor + r$ with $r < q$ and where $\lfloor . \rfloor$ denotes the integer part. From
this we may write:
$$
D_{M^{-j-1}x_{\bf 1}} f(x) = \sigma^{-q \lfloor \frac{j+1}{q} \rfloor} 
(D f)(x).M^{-r}x_{\bf 1}
$$   
and then 
$$
D_{M^{-j-1}x_{\bf 1}} f(x) \sim  \sigma^{-q\lfloor \frac{j+1}{q} 
\rfloor} 
(D f)(x).x_{r_j}
$$   
where $r_j$ depends on $j$. Making the same reasoning for any order $\mu$ of
differentiation and any direction $x_{\bf i}$, we get, in $L^p$:
$$
\|(\tilde D_{M^{-j-1}}^{\mu} f)(x)\|_{\Lp} \sim  \sigma^{-q|\mu|
\lfloor \frac{j+1}{q} \rfloor } 
\| (\tilde D^{\mu} f)(x)\|_{\Lp}.
$$   
We may thus conclude that
\begin{eqnarray*}
\|\tilde D^{\mu} (v_{j+1}(x)-v_j(x)) \|_{\Lp} &\sim& 
\|\tilde D_{M^{-j-1}}^{\mu} (v_{j+1}(x)-v_j(x)) \|_{\Lp}  
\sigma^{q|\mu|
(\times\lfloor \frac{j+1}{q} \rfloor)|}\\
&\lsim& 
\rho_{p,|\mu|}(S)^j m^{-\frac{j+1}{p}} \sigma^{q|\mu|
\lfloor \frac{j+1}{q} \rfloor} \|\tilde \Delta^{|\mu|} v^0\|_
 {(\lp)^{\tilde q_{|\mu|}}}. 
\end{eqnarray*}
To state the above result, we have used the fact that the joint spectral radius is
independent of the directions used for its computation.
Since we have the hypothesis that $\rho_{p,|\mu|}(S) 
\leq m^{\frac{1}{p}-\frac{s}{d}}$ for $s > |\mu|$, we get that   
\begin{eqnarray*}
\|\tilde D^{\mu} (v_{j+1}(x)-v_j(x)) \|_{\Lp} \lsim 
\sigma^{(|\mu|-s)j} \|\tilde \Delta^{|\mu|} v^0\|_
 {(\lp)^{\tilde q_{|\mu|}}},
\end{eqnarray*}
which tends to zero with $j$, and thus the limit function is in $W_N^p(\R^{d})$
$\blacksquare$.\\
A comparison between Theorem \ref{segaln} and \ref{segaln1} shows 
that when the subdivision scheme reproduces exactly polynomials, which is the
case of interpolatory subdivision schemes, the convergence is ensured provided 
$\Phi_0$ also exactly reproduces polynomials. When the subdivision scheme only 
reproduces polynomial the convergence is ensured provided that $\Phi_0$ is a 
box spline. Note also that the condition on the joint spectral radius is the 
same. We are currently investigating illustrative examples which involve the 
adaptation of the local averaging subdivision scheme proposed in \cite{Co} to
our non-separable context.  
\section{Conclusion}
We have addressed the issue of the definition of nonlinear
subdivision schemes associated to isotropic dilation matrix $M$. 
After the definition of the convergence concept of such operators, 
we have studied the convergence of these subdivision schemes in $L^p$ and in 
Sobolev spaces. 
Based on the study of the joint spectral radius of these operators, 
we have exhibited sufficient conditions for the convergence of the proposed 
subdivision schemes. This study has also brought into light the 
importance of an appropriate choice of $\Phi_0$ to define the limit 
function. In that context, box splines functions have shown to be a very 
interesting tool.

%%---------------------------------------------------------------------------%%

\end{document}